\documentclass[a4paper,12pt]{amsart}
\usepackage[a4paper,hscale=0.7,vscale=0.75,centering]{geometry}
\usepackage{fullpage}
\usepackage{url}

\newtheorem{theorem}{Theorem}[section]
\newtheorem{lemma}[theorem]{Lemma}

\theoremstyle{definition}

\newtheorem{example}[theorem]{Example}

\newtheorem*{proof1}{Proof of Lemma \ref{truncatingCoefficientsLemma}}
\newtheorem*{proof2}{Proof of Theorem \ref{globalOneSidedTheorem}}
\newtheorem*{proof3}{Proof of Lemma \ref{lemmaFeller}}
\newtheorem*{proof4}{Proof of Theorem \ref{existenceLOL}}
\newtheorem*{proof5}{Proof of Theorem \ref{existenceInvariant}}

\theoremstyle{remark}
\newtheorem{remark}[theorem]{Remark}

\numberwithin{equation}{section}

\begin{document}

\title{A note on existence of global solutions and invariant measures for jump SDEs with locally one-sided Lipschitz drift}

\author{Mateusz B. Majka}
\address{Institute for Applied Mathematics, University of Bonn, Endenicher Allee 60, 53115 Bonn, Germany}

\email{majka@uni-bonn.de}



\date{}

\dedicatory{}

\begin{abstract}
We extend some methods developed by Albeverio, Brze\'{z}niak and Wu and we show how to apply them in order to prove existence of global strong solutions of stochastic differential equations with jumps, under a local one-sided Lipschitz condition on the drift (also known as a monotonicity condition) and a local Lipschitz condition on the diffusion and jump coefficients, while an additional global one-sided linear growth assumption is satisfied. Then we use these methods to prove existence of invariant measures for a~broad class of such equations.
\end{abstract}

\maketitle

\section{Existence of global solutions under local Lipschitz conditions}

 Consider a stochastic differential equation in $\mathbb{R}^d$ of the form
 \begin{equation}\label{multiplicativeSDE}
  dX_t = b(X_t)dt + \sigma(X_t)dW_t + \int_{U} g(X_{t-},u) \widetilde{N}(dt,du) \,.
 \end{equation}
 Here $b : \mathbb{R}^d \to \mathbb{R}^d$, $\sigma : \mathbb{R}^d \to \mathbb{R}^{d \times d}$ and $g: \mathbb{R}^d \times U \to \mathbb{R}^d$, where $(W_t)_{t \geq 0}$ is a $d$-dimensional Wiener process, $(U, \mathcal{U}, \nu)$ is a $\sigma$-finite measure space and $N(dt,du)$ is a Poisson random measure on $\mathbb{R}^{+} \times U$ with intensity measure $dt \, \nu(du)$, while $\widetilde{N}(dt,du) = N(dt,du) - dt \, \nu(du)$ is the compensated Poisson random measure. We denote by $(X_t(x))_{t \geq 0}$ a solution to (\ref{multiplicativeSDE}) with initial condition $x \in \mathbb{R}^d$ and $\| \cdot \|_{HS}$ denotes the Hilbert-Schmidt norm of a matrix.
 The main result of the present paper is Theorem \ref{existenceInvariant}, where we prove existence of invariant measures for a certain class of such equations. However, we would first like to discuss the matter of existence of strong solutions to (\ref{multiplicativeSDE}), in the context of the paper \cite{albeverio} by Albeverio, Brze\'{z}niak and Wu.
 We claim that the following result holds.
 \begin{theorem}\label{existenceLOL}
  Assume that the coefficients in (\ref{multiplicativeSDE}) satisfy the following \emph{local one-sided Lipschitz} condition, i.e., for every $R > 0$ there exists $C_R > 0$ such that for any $x$, $y \in \mathbb{R}^d$ with $|x|$, $|y| \leq R$ we have
\begin{equation}\label{LOL}
 \langle b(x) - b(y), x - y \rangle + \| \sigma(x)  - \sigma(y) \|_{HS}^2 + \int_U |g(x,u) - g(y,u)|^2 \nu(du) \leq C_R |x-y|^2 \,.
\end{equation}
Moreover, assume a \emph{global one-sided linear growth} condition, i.e., there exists $C > 0$ such that for any $x \in \mathbb{R}^d$ we have
\begin{equation}\label{GOLG}
 \langle b(x) , x \rangle + \| \sigma(x) \|_{HS}^2 + \int_U |g(x,u)|^2 \nu(du) \leq C(1 + |x|^2) \,.
\end{equation}
Under (\ref{LOL}) and (\ref{GOLG}) and an additional assumption that $b: \mathbb{R}^d \to \mathbb{R}^d$ is continuous, there exists a unique global strong solution to (\ref{multiplicativeSDE}).
 \end{theorem}

 The one-sided Lipschitz condition (\ref{LOL}) above is sometimes called a monotonicity condition (see e.g. \cite{gyongykrylov} or \cite{krylovrozovski}) or a dissipativity condition (\cite{yma}, \cite{shao} or \cite{lwu}), although the term ``dissipativity'' is often reserved for the case in which (\ref{LOL}) is satisfied with a negative constant $C_R < 0$. We keep using the latter convention, calling (\ref{LOL}) one-sided Lipschitz regardless of the sign of the constant and using the term dissipativity only if the constant is negative. Note that the above theorem is a generalization of the following classic result.
 
 \begin{theorem}\label{existenceGL}
  Assume that the coefficients in (\ref{multiplicativeSDE}) satisfy a \emph{global Lipschitz} condition, i.e., there exists $C > 0$ such that for any $x$, $y \in \mathbb{R}^d$ we have
\begin{equation}\label{GL}
 | b(x) - b(y)|^2 + \| \sigma(x)  - \sigma(y) \|_{HS}^2 + \int_U |g(x,u) - g(y,u)|^2 \nu(du) \leq C |x-y|^2 \,.
\end{equation}
Moreover, assume a \emph{global linear growth} condition, i.e., there exists $L > 0$ such that for any $x \in \mathbb{R}^d$ we have
\begin{equation}\label{GLG}
 | b(x) |^2 + \| \sigma(x) \|_{HS}^2 + \int_U |g(x,u)|^2 \nu(du) \leq L(1 + |x|^2) \,.
\end{equation}
Under (\ref{GL}) and (\ref{GLG}) there exists a unique strong solution to (\ref{multiplicativeSDE}).
 \end{theorem}

 Theorem \ref{existenceGL} is very well-known and its proof can be found in many textbooks, see e.g. Theorem IV-9.1 in \cite{ikeda} or Theorem 6.2.3 in \cite{apple}. However, Theorem \ref{existenceLOL} is not so widespread in the literature and we had significant problems with finding a suitable reference for such a result. We finally learned that Theorem \ref{existenceLOL} can be inferred from Theorem 2 in \cite{gyongykrylov}, where a more general result is proved for equations driven by locally square integrable c\'{a}dl\'{a}g martingales taking values in Hilbert spaces.

Nevertheless, many authors use existence of solutions to equations like (\ref{multiplicativeSDE}) under a~one-sided Lipschitz condition for the drift (see e.g. \cite{yma}, \cite{shao}, \cite{jwang4}, \cite{lwu} for examples of some recent papers) claiming that this result is well-known, without giving any reference or while referring to positions that do not contain said result. Books that appear in this context include e.g. \cite{bichteler} and \cite{siturong} which, admittedly, contain various interesting extensions of the classic Theorem \ref{existenceGL}, but not the extension in which the Lipschitz condition is replaced with a one-sided Lipschitz condition and the linear growth with a one-sided linear growth.

Moreover, in a quite recent paper \cite{albeverio}, Albeverio, Brze\'{z}niak and Wu proved the following result (see Theorem 3.1 therein).

\begin{theorem}\label{albeverioTheorem}
 Assume that the coefficients in (\ref{multiplicativeSDE}) are such that for any $R > 0$ there exists $C_R > 0$ such that for any $x$, $y \in \mathbb{R}^d$ with $|x|$, $|y| \leq R$ we have
 \begin{equation}\label{albeverioLocalLipschitz}
  | b(x) - b(y)|^2 + \| \sigma(x)  - \sigma(y) \|_{HS}^2 \leq C_R |x-y|^2 \,.
 \end{equation}
Moreover, there exists $L > 0$ such that for any $x$, $y \in \mathbb{R}^d$ we have
\begin{equation}\label{albeverioGlobalLipschitz}
 \int_U |g(x,u) - g(y,u)|^2 \nu(du) \leq L |x-y|^2 \,.
\end{equation}
Finally, we assume a \emph{global one-sided linear growth} condition exactly like (\ref{GOLG}), i.e., there exists $C > 0$ such that for any $x \in \mathbb{R}^d$ we have
\begin{equation*}
 \langle b(x) , x \rangle + \| \sigma(x) \|_{HS}^2 + \int_U |g(x,u)|^2 \nu(du) \leq C(1 + |x|^2) \,.
\end{equation*}
Then there exists a unique global strong solution to (\ref{multiplicativeSDE}).
\end{theorem}

It is clear that Theorem \ref{albeverioTheorem} is less general than Theorem \ref{existenceLOL} and thus it is also a~special case of Theorem 2 in \cite{gyongykrylov}. Nevertheless, the proof in \cite{albeverio} is clearer and more direct than the one in \cite{gyongykrylov}, where the authors consider a much more general case. The main idea in \cite{albeverio} is to modify the locally Lipschitz coefficients in such a way as to obtain globally Lipschitz functions that agree with the given coefficients on a ball of fixed radius. Then using the classic Theorem \ref{existenceGL} it is possible to obtain a solution in every such ball and then to ``glue'' such local solutions by using the global one-sided linear growth condition to obtain a global solution. It is important to mention that the authors of \cite{albeverio} also use their methods to prove existence of invariant measures for a broad class of equations of the form (\ref{multiplicativeSDE}).

In view of all the above comments, we feel that it is necessary to give a direct proof of Theorem \ref{existenceLOL}. Following the spirit of the proof in \cite{albeverio}, we show how to extend the classic result (Theorem \ref{existenceGL}) in a step-by-step way in order to obtain Theorem \ref{existenceLOL}. Then we explain how to use the methods from \cite{albeverio} to obtain existence of invariant measures in our case, see Theorem \ref{existenceInvariant} in Section \ref{sectionInvariant}. The latter is an original result with potential applications in the theory of SPDEs, see Example \ref{exampleSPDE}.

For proving both Theorem \ref{existenceLOL} and \ref{existenceInvariant} we need the following auxiliary result regarding a possible modification of the coefficients in (\ref{multiplicativeSDE}).

\begin{lemma}\label{truncatingCoefficientsLemma}
 Assume that the coefficients in (\ref{multiplicativeSDE}) satisfy the local one-sided Lipschitz condition (\ref{LOL}) and that they are locally bounded in the sense that for every $R > 0$ there exists an $M_R > 0$ such that for all $x \in \mathbb{R}^d$ with $|x| \leq R$ we have
 \begin{equation}\label{localBoundedness}
  |b(x)|^2 + \| \sigma(x) \|_{HS}^2 + \int_U |g(x,u)|^2 \nu(du) \leq M_R \,.
 \end{equation}
Then for every $R > 0$ there exist truncated functions $b_R : \mathbb{R}^d \to \mathbb{R}^d$, $\sigma_R : \mathbb{R}^d \to \mathbb{R}^{d \times d}$ and $g_R : \mathbb{R}^d \times U \to \mathbb{R}^d$ such that for all $x \in \mathbb{R}^d$ with $|x| \leq R$ we have 
\begin{equation}\label{coeffAgree}
 b_R(x) = b(x) \text{, } \sigma_R(x) = \sigma(x) \text{ and } g_R(x,u) = g(x,u) \text{ for all } u \in \mathbb{R}^d \,.
\end{equation}
Moreover, $b_R$, $\sigma_R$ and $g_R$ satisfy a global one-sided Lipschitz condition, i.e., there exists a constant $C(R) > 0$ such that for all $x$, $y \in \mathbb{R}^d$ we have
\begin{equation}\label{GOLforTrunc}
  \langle b_R(x) - b_R(y), x - y \rangle + \| \sigma_R(x)  - \sigma_R(y) \|_{HS}^2 + \int_U |g_R(x,u) - g_R(y,u)|^2 \nu(du) \leq C(R) |x-y|^2
\end{equation}
and they are globally bounded, which means that there exists $M(R) > 0$ such that for all $x \in \mathbb{R}^d$ we have
\begin{equation}\label{coeffGloballyBounded}
 |b_R(x)|^2 + \| \sigma_R(x) \|_{HS}^2 + \int_U |g_R(x,u)|^2 \nu(du) \leq M(R) \,.
\end{equation}
\end{lemma}

Then, combining Theorem \ref{existenceGL} and Lemma \ref{truncatingCoefficientsLemma}, we are able to prove existence of solutions while the coefficients in (\ref{multiplicativeSDE}) are bounded and satisfy a global one-sided Lipschitz condition.

\begin{theorem}\label{globalOneSidedTheorem}
 Assume that $b$ is continuous and that the coefficients in (\ref{multiplicativeSDE}) satisfy a~global one-sided Lipschitz condition, i.e., there exists $K > 0$ such that for all $x$, $y \in \mathbb{R}^d$ we have
 \begin{equation}\label{GOL}
  \langle b(x) - b(y), x - y \rangle + \| \sigma(x)  - \sigma(y) \|_{HS}^2 + \int_U |g(x,u) - g(y,u)|^2 \nu(du) \leq K |x-y|^2 \,.
 \end{equation}
Additionally, assume that the coefficients are globally bounded, i.e., there exists $M > 0$ such that for all $x \in \mathbb{R}^d$ we have
\begin{equation}\label{boundedCoefficients}
 |b(x)|^2 + \| \sigma(x) \|_{HS}^2 + \int_U |g(x,u)|^2 \nu(du) \leq M \,.
\end{equation}
Then there exists a unique strong solution to (\ref{multiplicativeSDE}).
\end{theorem}

The proofs of Lemma \ref{truncatingCoefficientsLemma} and Theorem \ref{globalOneSidedTheorem} can be found in Section \ref{sectionProofs}. Having proved the above two results, we proceed with the proof of Theorem \ref{existenceLOL} as in \cite{albeverio} (see Proposition 2.9 and Theorem 3.1 therein, see also \cite{gyongykrylov}, page 14, for a similar reasoning). More details can be found at the end of Section \ref{sectionProofs} below.

\section{Existence of invariant measures}\label{sectionInvariant}

 The existence of an invariant measure for the solution of (\ref{multiplicativeSDE}) is shown using the Krylov-Bogoliubov method, see e.g. Theorem III-2.1 in \cite{hasminskii} and the discussion in the introduction to \cite{dapratogatarek}. It follows from there that for the existence of an invariant measure for a process $(X_t)_{t \geq 0}$ with a Feller semigroup $(p_t)_{t \geq 0}$ it is sufficient to show that for some $x \in \mathbb{R}^d$ the process $(X_t(x))_{t \geq 0}$ is bounded in probability at infinity in the sense that for any $\varepsilon > 0$ there exist $R > 0$ and $t > 0$ such that for all $s \geq t$ we have
 \begin{equation}\label{boundedInProbability}
  \mathbb{P}(|X_s(x)| > R) < \varepsilon \,.
 \end{equation}
Therefore if we show that there exist constants $M$, $K > 0$ such that
\begin{equation}\label{theorem45ineq}
 \mathbb{E}|X_t(x)|^2 \leq |x|^2 e^{-Kt} + M/K
\end{equation}
holds for all $t \geq 0$, then (\ref{boundedInProbability}) follows easily by the Chebyshev inequality and we obtain the existence of an invariant measure. Based on this idea, we can prove the following result.

\begin{theorem}\label{existenceInvariant}
 Assume that the coefficients in (\ref{multiplicativeSDE}) satisfy the local one-sided Lipschitz condition (\ref{LOL}) and that there exist constants $K$, $M > 0$ such that for all $x \in \mathbb{R}^d$ we have
 \begin{equation}\label{dissipativeGrowth}
   \langle b(x) , x \rangle + \| \sigma(x) \|_{HS}^2 + \int_U |g(x,u)|^2 \nu(du) \leq -K|x|^2 + M \,.
 \end{equation}
Assume also that there exists a constant $L > 0$ such that for all $x \in \mathbb{R}^d$ we have
\begin{equation}\label{linearGrowthSeparate}
 \| \sigma(x) \|_{HS}^2 + \int_U |g(x,u)|^2 \nu(du) \leq L(1 + |x|^2) \,.
\end{equation}
Finally, let the drift coefficient $b$ in (\ref{multiplicativeSDE}) be continuous. Then there exists an invariant measure for the solution of (\ref{multiplicativeSDE}).
\end{theorem}

We can compare this result with the one proved in \cite{albeverio} (see Theorem 4.5 therein).

\begin{theorem}
 Assume that the coefficients $b$ and $\sigma$ in (\ref{multiplicativeSDE}) satisfy the local Lipschitz condition (\ref{albeverioLocalLipschitz}) and that $g$ satisfies the global Lipschitz condition (\ref{albeverioGlobalLipschitz}). Assume also the condition (\ref{dissipativeGrowth}) as in the Theorem \ref{existenceInvariant} above. Then there exists an invariant measure for the solution of (\ref{multiplicativeSDE}).
\end{theorem}

\begin{remark}\label{remarkSeparateLinearGrowth}
 Observe that our additional condition (\ref{linearGrowthSeparate}) in Theorem \ref{existenceInvariant} does not follow from (\ref{dissipativeGrowth}) since $\langle b(x) , x \rangle$ can be negative. Therefore it would seem that our result is not a straightforward generalization of Theorem 4.5 in \cite{albeverio}. However, we believe that the condition (\ref{linearGrowthSeparate}) is also necessary to prove Theorem 4.5 in \cite{albeverio}, at least we were not able to retrace the proof of Proposition 4.3 therein (which is crucial for the proof of Theorem 4.5) without this additional condition. Therefore we are convinced that (\ref{linearGrowthSeparate}) should be added to the list of assumptions of Theorem 4.5 in \cite{albeverio} and that our result is indeed its strict generalization. This has been confirmed in our private communication with one of the authors of \cite{albeverio}.
\end{remark}

For the proof of Theorem \ref{existenceInvariant} we first need the following fact, which can be proved exactly like in \cite{albeverio}.

\begin{lemma}
 The solution $(X_t)_{t \geq 0}$ to the equation (\ref{multiplicativeSDE}) is a strong Markov process and thus it generates a Markov semigroup $(p_t)_{t \geq 0}$.
 \begin{proof}
  See Proposition 4.1 and Proposition 4.2 in \cite{albeverio}.
 \end{proof}
\end{lemma}

Now we need the following lemma, which is a generalization of Proposition 4.3 in \cite{albeverio} (see Remark \ref{remarkSeparateLinearGrowth} about inclusion of the assumption (\ref{linearGrowthSeparate})).

\begin{lemma}\label{lemmaFeller}
 Under the assumptions (\ref{LOL}), (\ref{GOLG}), (\ref{linearGrowthSeparate}) and if $b$ is continuous, the semigroup $(p_t)_{t \geq 0}$ associated with the solution $(X_t)_{t \geq 0}$ of (\ref{multiplicativeSDE}) is Feller.
\end{lemma}

Having proved the above lemma, we can easily conclude the proof of Theorem \ref{existenceInvariant}, following the proof of Theorem 4.5 in \cite{albeverio}, i.e., we just use the condition (\ref{dissipativeGrowth}) to show (\ref{theorem45ineq}) and then use the Krylov-Bogoliubov method presented above. More details can be found in Section \ref{sectionProofs}.

Before concluding this section, let us look at some examples.

\begin{example}
 Consider an SDE of the form (\ref{multiplicativeSDE}) with the drift given by
 \begin{equation*}
  b(x):= -x |x|^{-\alpha} \mathbf{1}_{\{ x \neq 0 \} } \,,
 \end{equation*}
where $\alpha \in (0,1)$. Equations of this type are considered in Example 171 in \cite{siturong}. It is easy to check that the function $b$ defined above is not locally Lipschitz, since it does not satisfy a Lipschitz condition in any neighbourhood of zero. However, we can show that it satisfies a one-sided Lipschitz condition globally with constant zero. Indeed, following the calculations in Example 171 in \cite{siturong}, for any nonzero $x$, $y \in \mathbb{R}^d$ we have
\begin{equation*}
 \begin{split}
  \langle x - y , -x|x|^{-\alpha} + y|y|^{-\alpha} \rangle &= -|x|^{2 - \alpha} + \langle y , x|x|^{-\alpha} \rangle +  \langle x , y|y|^{-\alpha} \rangle -|y|^{2 - \alpha} \\
  &\leq -|x|^{2 - \alpha} -|y|^{2 - \alpha} + |y||x|^{1 - \alpha} + |x||y|^{1 - \alpha} \\
  &= (|x| - |y|)(|y|^{1 - \alpha} - |x|^{1 - \alpha}) \leq 0 \,,
 \end{split}
\end{equation*}
where the last inequality holds since $1 - \alpha \in (0,1)$. Thus, if we consider an equation of the form (\ref{multiplicativeSDE}) with the drift $b$ and any locally Lipschitz coefficients $\sigma$ and $g$, the condition (\ref{LOL}) is satisfied. Moreover, if $\sigma$ and $g$ satisfy the global linear growth condition (\ref{linearGrowthSeparate}) with some constant $L > 0$, then by replacing the drift $b$ defined above with
\begin{equation*}
 \widetilde{b}(x) := b(x) - Kx \,,
\end{equation*}
where $K > L$, we obtain coefficients that satisfy (\ref{dissipativeGrowth}). More generally, we can take
\begin{equation*}
  \widetilde{b}(x) := b(x) - \nabla U(x) \,,
\end{equation*}
where $U$ is a strongly convex function with convexity constant $K > L$. This way we obtain a class of examples of equations for which our Theorem \ref{existenceInvariant} applies, but Theorem 4.5 in \cite{albeverio} does not, since the local Lipschitz assumption is not satisfied.
\end{example}

\begin{example}\label{exampleSPDE}
 Our results may have applications in the study of stochastic evolution equations with L\'{e}vy noise on infinite dimensional spaces, where the coefficients are often not Lipschitz, see e.g. \cite{brzezniak} and the references therein. In particular, in \cite{brzezniak} the authors consider SPDEs with drifts satisfying a local monotonicity condition and use their finite dimensional approximations, which may lead to SDEs satisfying our condition (\ref{LOL}), cf. the condition (H2) and the formula (4.4) in \cite{brzezniak}.
\end{example}

\section{Proofs}\label{sectionProofs}

In order to keep our presentation compact, we will only present the proof of Theorem \ref{existenceLOL} in a slightly less general setting than that presented in the first section. Namely, we will additionally assume that the diffusion coefficient $\sigma$ and the jump coefficient $g$ in the equation (\ref{multiplicativeSDE}) satisfy a local Lipschitz condition separately from the drift $b$, i.e., for every $R > 0$ there exists $S_R > 0$ such that for any $x$, $y \in \mathbb{R}^d$ with $|x|$, $|y| \leq R$ we have
\begin{equation}\label{SeparateLL}
 \| \sigma(x)  - \sigma(y) \|_{HS}^2 + \int_U |g(x,u) - g(y,u)|^2 \nu(du) \leq S_R |x-y|^2 \,.
\end{equation}
Obviously, (\ref{SeparateLL}) does not follow from (\ref{LOL}), since the values of $\langle b(x) - b(y), x - y \rangle$ can be negative. However, requiring the condition (\ref{SeparateLL}) to be satisfied seems to be rather natural in many cases. It is possible to weaken this assumption and prove the exact statement of Theorem \ref{existenceLOL} using methods from Section 3 of Chapter II in \cite{krylovrozovski} (see also Section 3 in \cite{gyongykrylov}), but this creates additional technical difficulties and thus we decided to omit this extension here, aiming at a clear and straightforward presentation.

The consequence of adding the assumption (\ref{SeparateLL}) is that the coefficients of (\ref{multiplicativeSDE}) automatically satisfy the local boundedness condition (\ref{localBoundedness}) required in Lemma \ref{truncatingCoefficientsLemma} (remember that $b$ is assumed to be continuous and thus it is locally bounded anyway). It also means that from Lemma \ref{truncatingCoefficientsLemma} we obtain coefficients $\sigma_R$ and $g_R$ that satisfy a separate global Lipschitz condition, i.e., the condition (\ref{GOLforTrunc}) without the term involving $b_R$. Hence we can prove Theorem \ref{globalOneSidedTheorem} under an additional assumption, i.e., we can use the fact that there exists $S > 0$ such that for all $x$, $y \in \mathbb{R}^d$ we have
\begin{equation}\label{SeparateGL}
  \| \sigma(x)  - \sigma(y) \|_{HS}^2 + \int_U |g(x,u) - g(y,u)|^2 \nu(du) \leq S |x-y|^2 \,.
\end{equation}

However, the assumption (\ref{SeparateLL}) is not needed for the proof of Theorem \ref{existenceInvariant}, where we also use Lemma \ref{truncatingCoefficientsLemma}, but we do not need to obtain truncated coefficients $\sigma_R$ and $g_R$ satisfying a separate global Lipschitz condition and the assumption about local boundedness is guaranteed by the separate linear growth condition (\ref{linearGrowthSeparate}) and the continuity of $b$. Thus the reasoning presented below gives a complete proof of the exact statement of our Theorem \ref{existenceInvariant}.

\begin{proof1}
 For a related reasoning, see the proof of Lemma 4 in \cite{gyongykrylov} or Lemma 172 in \cite{siturong}. Note that the method of truncating the coefficients of (\ref{multiplicativeSDE}) which was used in the proof of Proposition 2.7 in \cite{albeverio} and which works in the case of Lipschitz coefficients, does not work for a one-sided Lipschitz drift and thus we need a different approach. For any $R > 0$, we can consider a smooth, non-negative function with compact support $\eta_R \in \mathcal{C}^{\infty}_{c}(\mathbb{R}^d)$ such that
 \begin{equation*}
  \eta_R(x) = \begin{cases}
                 1 \,, & \text{if } |x| \leq R \,, \\
                 0 \,, & \text{if } |x| > R + 1 
                \end{cases}
 \end{equation*}
and $\eta_R(x) \leq 1$ for all $x \in \mathbb{R}^d$. Then we can define 
\begin{equation*}
b_R(x) := \eta_R(x) b(x) \text{, } \sigma_R(x) := \eta_R(x) \sigma(x) \text{ and } g_R(x,u) := \eta_R(x) g(x,u) \text{ for all } u \in \mathbb{R}^d \,.
\end{equation*}
Then it is obvious that the condition (\ref{coeffAgree}) is satisfied and the condition (\ref{coeffGloballyBounded}) immediately follows from (\ref{localBoundedness}). Therefore it remains to be shown that the functions $b_R$, $\sigma_R$ and $g_R$ satisfy the global one-sided Lipschitz condition (\ref{GOLforTrunc}). We have
\begin{equation}\label{lemmaCalculations}
\begin{split}
   \langle b_R(x) &- b_R(y), x - y \rangle + \| \sigma_R(x)  - \sigma_R(y) \|_{HS}^2 + \int_U |g_R(x,u) - g_R(y,u)|^2 \nu(du) \\
   &= \langle \eta_R(x)b(x) - \eta_R(y)b(y), x - y \rangle + \| \eta_R(x)\sigma(x)  - \eta_R(y)\sigma(y) \|_{HS}^2 \\
   &+ \int_U |\eta_R(x)g(x,u) - \eta_R(y)g(y,u)|^2 \nu(du) \\
   &\leq \eta_R(x)\langle b(x) - b(y), x - y \rangle + \langle (\eta_R(x) - \eta_R(y)) b(y), x - y \rangle \\
   &+ |\eta_R(x)|^2 \| \sigma(x)  - \sigma(y) \|_{HS}^2 + \| (\eta_R(x) - \eta_R(y)) \sigma(y) \|_{HS}^2 \\
   &+ \int_U |\eta_R(x)|^2 |g(x,u) - g(y,u)|^2 \nu (du) + \int_U |\eta_R(x) - \eta_R(y)|^2 |g(y,u)|^2 \nu (du) \,.
\end{split}
\end{equation}
Now assume $x$ and $y$ are such that
\begin{equation}\label{etaAssumption}
 \eta_R(y) \geq \eta_R(x) > 0 \,.
\end{equation}
The case when $\eta_R(x) = 0$ is simpler and the case $\eta_R(y) \leq \eta_R(x)$ can be handled by changing the role of $x$ and $y$ in the calculations above. From (\ref{etaAssumption}) it follows that $|y| \leq R + 1$ and $|x| \leq R + 1$ and thus we can use the local one-sided Lipschitz condition (\ref{LOL}) with $R + 1$ to get
\begin{equation*}
 \langle b(x) - b(y), x - y \rangle + \| \sigma(x)  - \sigma(y) \|_{HS}^2 + \int_U |g(x,u) - g(y,u)|^2 \nu(du) \leq C_{R+1} |x-y|^2 
\end{equation*}
with some constant $C_{R+1}$. Combining this with the fact that $\eta_R \leq 1$ (and thus $\eta_R^2 \leq \eta_R$) allows us to bound the sum of the first, the third and the fifth term on the right hand side of (\ref{lemmaCalculations}) by $C_{R+1}|x-y|^2$. Observe now that the function $\eta_R$ is Lipschitz (with a constant, say, $C_{Lip(\eta_R)}$) and thus
\begin{equation*}
 \eta_R(x) - \eta_R(y) \leq C_{Lip(\eta_R)} |x-y| \,.
\end{equation*}
Since $|y| \leq R + 1$, we can use the local boundedness condition (\ref{localBoundedness}) with some constant $M_{R+1}$. We first bound $|b(y)|$ by the square root of the left hand side of (\ref{localBoundedness}) in order to get
\begin{equation*}
 \langle (\eta_R(x) - \eta_R(y)) b(y), x - y \rangle \leq \sqrt{M_{R+1}} C_{Lip(\eta_R)} |x-y|^2 \,.
\end{equation*}
Then we use (\ref{localBoundedness}) once again in order to bound the sum of the fourth and the sixth term on the right hand side of (\ref{lemmaCalculations}) by $M_{R+1} C_{Lip(\eta_R)}^2 |x-y|^2$. Combining all these facts together, we can bound the right hand side of (\ref{lemmaCalculations}) by 
\begin{equation*}
 C_{R+1}|x-y|^2 + \sqrt{M_{R+1}} C_{Lip(\eta_R)} |x-y|^2 +  M_{R+1} C_{Lip(\eta_R)}^2 |x-y|^2 \,.
\end{equation*}
Therefore the global one-sided Lipschitz condition for $b_R$, $\sigma_R$ and $g_R$ is satisfied with a~constant 
\begin{equation*}
 C(R) := C_{R+1} + \sqrt{M_{R+1}} C_{Lip(\eta_R)} +  M_{R+1} C_{Lip(\eta_R)}^2 \,,
\end{equation*}
which finishes the proof.
\qed
\end{proof1}

Before we proceed to the proof of Theorem \ref{globalOneSidedTheorem}, let us formulate a crucial technical lemma. Its proof is just a slightly altered second part of the proof of Lemma 3.3 in \cite{yma}, but we include the full calculations here for completeness and, more importantly, because we need to use a related, but modified reasoning in the proof of Theorem \ref{globalOneSidedTheorem}. The lemma itself will be used in the proof of Lemma \ref{lemmaFeller} later on.

\begin{lemma}\label{ymaEstimates}
 Assume that the coefficients of the equation (\ref{multiplicativeSDE}) with an initial condition $x \in \mathbb{R}^d$ satisfy the global one-sided linear growth condition (\ref{GOLG}) and that $\sigma$ and $g$ additionally satisfy the separate linear growth condition (\ref{linearGrowthSeparate}). Then there exist constants $\widetilde{C} > 0$ and $\widetilde{K} > 0$ such that
 \begin{equation*}
  \mathbb{E} \sup_{s \leq t} |X_s|^2 \leq \widetilde{K} e^{2\widetilde{C}t}(1+|x|^2) \,,
 \end{equation*}
 where $(X_t)_{t \geq 0} = (X_t(x))_{t \geq 0}$ is a solution to (\ref{multiplicativeSDE}) with initial condition $x \in \mathbb{R}^d$.
\begin{proof}
 By the It\^{o} formula, we have
 \begin{equation}\label{ItoYma}
  \begin{split}
   |X_t|^2 &= |x|^2 + 2\int_0^t \langle b(X_s), X_s \rangle ds + 2\int_0^t \langle \sigma(X_s), X_s dW_s \rangle \\
   &+ \int_0^t \| \sigma(X_s) \|_{HS}^2 ds +2 \int_0^t \int_U \langle g(X_{s-},u) , X_s \rangle \widetilde{N}(ds,du) \\
   &+ \int_0^t \int_U  |g(X_{s-},u)|^2 N(ds,du) \,.
  \end{split}
 \end{equation}
Now let us consider the process
\begin{equation*}
 M_t := \int_0^t \langle \sigma(X_s), X_s dW_s \rangle + \int_0^t \int_U \langle g(X_{s-},u) , X_s \rangle \widetilde{N}(ds,du) \,,
\end{equation*}
which is a local martingale. Thus, by the Burkholder-Davis-Gundy inequality, there exists a constant $C_1 > 0$ such that
\begin{equation*}
 \begin{split}
  \mathbb{E} \sup_{s \leq t} |M_s| &\leq C_1 \mathbb{E} \left[ \int_0^t |\sigma^{*}(X_s)X_s|^2 ds + \int_0^t \int_U |\langle g(X_{s-},u) , X_s \rangle|^2 N(ds,du)
  \right]^{\frac{1}{2}} \\
  &\leq C_1 \mathbb{E} \left[ (\sup_{s \leq t} |X_s|^2 ) \left( \int_0^t \| \sigma^{*}(X_s) \|^2 ds + \int_0^t \int_U | g(X_{s-},u)|^2 N(ds,du) \right)
  \right]^{\frac{1}{2}} \\
  &\leq C_1 \left(\mathbb{E} \sup_{s \leq t} |X_s|^2  \right)^{\frac{1}{2}} \left(\mathbb{E} \left[ \int_0^t \| \sigma^{*}(X_s) \|^2 ds + \int_0^t \int_U | g(X_{s-},u)|^2 N(ds,du) \right] \right)^{\frac{1}{2}} \\
  &\leq \frac{C_1}{2}a \mathbb{E} \sup_{s \leq t} |X_s|^2 + \frac{C_1}{2a}\mathbb{E} \left[ \int_0^t \| \sigma^{*}(X_s) \|^2 ds + \int_0^t \int_U | g(X_{s-},u)|^2 N(ds,du) \right] \\
  &\leq \frac{C_1}{2}a \mathbb{E} \sup_{s \leq t} |X_s|^2 + \frac{C_1}{2a} L \mathbb{E} \int_0^t (|X_s|^2 + 1)ds \,.
 \end{split}
\end{equation*}
Here $\| \cdot \|$ denotes the operator norm and $\sigma^{*}$ is a transposed $\sigma$. In the third step we used the H\"{o}lder inequality in the form $\mathbb{E}A^{\frac{1}{2}}B^{\frac{1}{2}} \leq (\mathbb{E}A)^{\frac{1}{2}}(\mathbb{E}B)^{\frac{1}{2}}$, in the fourth step we used $(AB)^{\frac{1}{2}} \leq \frac{1}{2}aA + \frac{1}{2a}B$ for any $a > 0$, which can be chosen later, and in the fifth step we used the separate global linear growth condition (\ref{linearGrowthSeparate}) for $\sigma$ and $g$ along with the fact that $\| \cdot \| \leq \| \cdot \|_{HS}$. Now we can use the formula (\ref{ItoYma}) to get
\begin{equation}\label{ymaMain}
 \begin{split}
   \mathbb{E} \sup_{s \leq t} |X_s|^2 &\leq |x|^2 + 2\mathbb{E} \sup_{s \leq t} |M_s| + 2\mathbb{E} \sup_{s \leq t} \int_0^s \langle b(X_r), X_r \rangle dr \\
   &+ \mathbb{E} \sup_{s \leq t} \left[ \int_0^s \| \sigma(X_r) \|_{HS}^2 dr + \int_0^s \int_U  |g(X_{r-},u)|^2 N(dr,du) \right] \,.
 \end{split}
\end{equation}
Observe that obviously 
\begin{equation}\label{separateOLforB}
 \langle b(X_r), X_r \rangle \leq \langle b(X_r), X_r \rangle + \| \sigma(X_r) \|_{HS}^2 + \int_U  |g(X_{r-},u)|^2 \nu(du)
\end{equation}
and thus from the global one-sided linear growth condition (\ref{GOLG}) we get
\begin{equation*}
 \mathbb{E} \sup_{s \leq t} \int_0^s \langle b(X_r), X_r \rangle dr \leq C \mathbb{E} \sup_{s \leq t} \int_0^s (|X_r|^2 + 1) dr \leq C \mathbb{E} \int_0^t (|X_r|^2 + 1) dr \,.
\end{equation*}
On the other hand, using the separate linear growth condition (\ref{linearGrowthSeparate}) we get
\begin{equation*}
\begin{split}
 \mathbb{E} \sup_{s \leq t} &\left[ \int_0^s \| \sigma(X_r) \|_{HS}^2 dr + \int_0^s \int_U  |g(X_{r-},u)|^2 N(dr,du) \right] \\
 &= \mathbb{E} \left[ \int_0^t \| \sigma(X_r) \|_{HS}^2 dr + \int_0^t \int_U  |g(X_{r-},u)|^2 N(dr,du) \right] \\
 &= \mathbb{E} \left[ \int_0^t \| \sigma(X_r) \|_{HS}^2 dr + \int_0^t \int_U  |g(X_{r-},u)|^2 \nu(du) dr \right] \\
 &\leq L \mathbb{E} \int_0^t (|X_r|^2 + 1) dr \,.
\end{split}
\end{equation*}
Combining all the above estimates, we get from (\ref{ymaMain}) that
\begin{equation*}
  \mathbb{E} \sup_{s \leq t} |X_s|^2 \leq |x|^2 + C_1 a \mathbb{E} \sup_{s \leq t} |X_s|^2 + (\frac{C_1}{a} L  + 2C +L)\mathbb{E} \int_0^t (|X_r|^2 + 1)dr \,.
\end{equation*}
Now, choosing $a = 1 / (2C_1)$ we obtain
\begin{equation*}
 \mathbb{E} \sup_{s \leq t} |X_s|^2 \leq 2 |x|^2 + 2(2C_1^2 L  + 2C +L)\mathbb{E} \int_0^t \sup_{w \leq r} (|X_w|^2 + 1)dr \,.
\end{equation*}
Hence, using the Gronwall inequality for the function $\mathbb{E} \sup_{s \leq t} |X_s|^2 + 1$ we get
\begin{equation*}
 \mathbb{E} \sup_{s \leq t} |X_s|^2 + 1 \leq 2 (|x|^2 + 1)  \exp (2(2C_1^2 L  + 2C +L)t) \,,
\end{equation*}
which finishes the proof.
\end{proof}

\end{lemma}

\begin{proof2}
 Let $j \in \mathcal{C}^{\infty}_{c}(\mathbb{R}^d)$ be a smooth function with a compact support contained in $B(0,1)$, such that $\int_{\mathbb{R}^d} j(z) dz = 1$. Then, for any $k \geq 1$, define
 \begin{equation*}
  b^k(x) := \int_{\mathbb{R}^d} b(x - \frac{z}{k}) j(z) dz \,.
 \end{equation*}
Now we can consider a sequence of equations
\begin{equation}\label{eqXk}
 dX^k_t = b^k(X^k_t)dt + \sigma(X^k_t)dW_t + \int_{U} g(X^k_{t-},u) \widetilde{N}(dt,du) \,.
\end{equation}
Note that we have replaced only the drift coefficient $b$ with $b^k$ while $\sigma$ and $g$ remain unchanged. This is due to the fact that we decided to prove Theorem \ref{existenceLOL} with an additional assumption of separate local Lipschitz condition (\ref{SeparateLL}) for $\sigma$ and $g$. Thanks to this, we can work in the present proof under an additional assumption that $\sigma$ and $g$ are globally Lipschitz, i.e. they satisfy (\ref{SeparateGL}), cf. the discussion at the beginning of this section. Now observe that the function $b^k$ defined above is also globally Lipschitz. Indeed, for any $x$, $y \in \mathbb{R}^d$ we have
\begin{equation*}
\begin{split}
|b^k(x) - b^k(y)| &= |\int_{\mathbb{R}^d} b(x - \frac{z}{k}) j(z) dz - \int_{\mathbb{R}^d} b(y - \frac{z}{k}) j(z) dz| \\
&= |k^d \int_{\mathbb{R}^d} b(w) j(k(x-w)) dz - k^d \int_{\mathbb{R}^d} b(w) j(k(y-w)) dw| \\
&\leq k^d \int_{\mathbb{R}^d} |b(w)| |j(k(x-w)) - j(k(y-w))| dw \\
&\leq k^{d+1} \sqrt{M} |x-y| \int_{\mathbb{R}^d} \sup_{w \in \mathbb{R}^d}|\nabla j(w)|dw \,,
\end{split}
\end{equation*}
where in the last step we use the fact that $b$ is bounded by $\sqrt{M}$ (cf. (\ref{boundedCoefficients})) and $j$~is Lipschitz with the Lipschitz constant given by the supremum of the norm of its gradient (which is obviously integrable since $j \in \mathcal{C}^{\infty}_{c}(\mathbb{R}^d)$). Having proved that $b^k$ is globally Lipschitz, we can use Theorem \ref{existenceGL} to ensure existence of a unique strong solution $(X_t^k)_{t \geq 0}$ to the equation (\ref{eqXk}). We will prove now that the sequence of solutions $\{(X_t^k)_{t \geq 0}\}_{k=1}^{\infty}$ has a limit (in the sense of almost sure convergence, uniform on bounded time intervals) and that this limit is in fact a solution to (\ref{multiplicativeSDE}). To this end, we will make use of the calculations from Lemma \ref{ymaEstimates}. 

Observe that for any $k$, $l \geq 1$, if we use the It\^{o} formula to calculate $|X_t^k - X_t^l|^2$, we will obtain exactly the formula (\ref{ItoYma}) with $X_t$ replaced by the difference $X_t^k - X_t^l$ and the function $b(X_s)$ replaced by $b^k(X_s^k) - b^l(X_s^l)$. Furthermore, we can make the term $|x|^2$ vanish (we can assume that all the solutions $(X_t^k)_{t \geq 0}$ have the same initial condition). Now we can proceed exactly like in the proof of Lemma \ref{ymaEstimates}, this time using the separate global Lipschitz condition (\ref{SeparateGL}) for $\sigma$ and $g$ in the steps where we used the separate linear growth condition (\ref{linearGrowthSeparate}) before, in order to get
\begin{equation}\label{proof2Estimate}
 \begin{split}
  \mathbb{E} \sup_{s \leq t} |X_s^k - X_s^l|^2 &\leq C_1 a \mathbb{E} \sup_{s \leq t} |X_s^k - X_s^l|^2 + (\frac{C_1}{a} S + S)\mathbb{E} \int_0^t |X_r^k - X_r^l|^2 dr \\
  &+ 2 \mathbb{E} \sup_{s \leq t} \int_0^s \langle X_r^k - X_r^l , b^k(X_r^k) - b^l(X_r^l) \rangle dr \,.
 \end{split}
\end{equation}
Thus the only term, with which we have to deal in a different way compared to the proof of Lemma \ref{ymaEstimates}, is the last one. We have
\begin{equation*}
 \begin{split}
  \mathbb{E} &\sup_{s \leq t} \int_0^s \langle X_r^k - X_r^l , b^k(X_r^k) - b^l(X_r^l) \rangle dr \\
  &= \mathbb{E} \sup_{s \leq t} \int_0^s \left\langle X_r^k - X_r^l , \int_{\mathbb{R}^d} b(X_r^k - \frac{z}{k})j(z)dz - \int_{\mathbb{R}^d} b(X_r^l - \frac{z}{l})j(z)dz\right\rangle dr \\
  &= \mathbb{E} \sup_{s \leq t} \int_0^s \Bigg\{  \int_{\mathbb{R}^d} \left\langle (X_r^k - \frac{z}{k}) - (X_r^l - \frac{z}{l}) , b(X_r^k - \frac{z}{k}) - b(X_r^l - \frac{z}{l}) \right\rangle j(z)dz  \\
  &+ \int_{\mathbb{R}^d} \left\langle \frac{z}{k} - \frac{z}{l} , b(X_r^k - \frac{z}{k}) - b(X_r^l - \frac{z}{l}) \right\rangle j(z)dz \Bigg\} dr \\
  &=: \mathbb{E} \sup_{s \leq t} \int_0^s (I_r^1 + I_r^2) dr \,.
 \end{split}
\end{equation*}
Now observe that since $b$ is assumed to be bounded by $\sqrt{M}$ (see (\ref{boundedCoefficients})), we have
\begin{equation*}
\begin{split}
 I_r^2 &\leq 2\sqrt{M} \int_{\mathbb{R}^d} \left|\frac{z}{k} - \frac{z}{l}\right| j(z) dz = 2\sqrt{M} \left|\frac{1}{k} - \frac{1}{l}\right| \int_{\mathbb{R}^d} |z|j(z) dz \\
 &=: 2\sqrt{M} \left|\frac{1}{k} - \frac{1}{l}\right| C^1(j) < \infty \,.
 \end{split}
\end{equation*}
As for $I_r^1$, we can use the one-sided Lipschitz condition (\ref{GOL}) for $b$ similarly like we used one-sided linear growth in (\ref{separateOLforB}) to get
\begin{equation*}
 \begin{split}
  I_r^1 &\leq K \int_{\mathbb{R}^d} \left|(X_r^k - \frac{z}{k}) - (X_r^l - \frac{z}{l})\right|^2 j(z)dz \\
  &= K \int_{\mathbb{R}^d} \left|(X_r^k - X_r^l) - (\frac{1}{k} - \frac{1}{l})z\right|^2 j(z)dz \\
  &\leq 2K \int_{\mathbb{R}^d} \left|X_r^k - X_r^l\right|^2 j(z)dz + 2K \int_{\mathbb{R}^d} \left|\frac{1}{k} - \frac{1}{l}\right|^2 |z|^2 j(z)dz \\
  &=: 2K \left|X_r^k - X_r^l\right|^2 + 2K C^2(j) \left|\frac{1}{k} - \frac{1}{l}\right|^2 \,.
 \end{split}
\end{equation*}
Combining the above estimates, we have
\begin{equation*}
\begin{split}
 \mathbb{E} \sup_{s \leq t} \int_0^s (I_r^1 + I_r^2) dr &\leq 2K \mathbb{E} \int_0^t \left|X_r^k - X_r^l\right|^2 dr + 2 t K C^2(j) \left|\frac{1}{k} - \frac{1}{l}\right|^2 \\
 &+ 2 t \sqrt{M} \left|\frac{1}{k} - \frac{1}{l}\right| C^1(j) \\
 &\leq 2 K \mathbb{E} \int_0^t \sup_{w \leq r} |X_w^k - X_w^l|^2 dr + \widehat{C} t \left|\frac{1}{k} - \frac{1}{l}\right| \,,
 \end{split}
\end{equation*}
where the last inequality holds with a constant $\widehat{C} := 2 K C^2(j) + 2 \sqrt{M} C^1(j)$ for $k$ and $l$ large enough so that $|\frac{1}{k} - \frac{1}{l}| < 1$. Now we can come back to (\ref{proof2Estimate}) and, taking $a = 1 / (2C_1)$, similarly like in the proof of Lemma \ref{ymaEstimates} we obtain
\begin{equation*}
 \begin{split}
  \mathbb{E} \sup_{s \leq t} |X_s^k - X_s^l|^2 &\leq 2(2C_1^2 S + S)\mathbb{E} \int_0^t \sup_{w \leq r} |X_w^k - X_w^l|^2 dr \\
  &+ 8 K \mathbb{E} \int_0^t \sup_{w \leq r} |X_w^k - X_w^l|^2 dr + 4 \widehat{C} t \left|\frac{1}{k} - \frac{1}{l}\right| \,.
 \end{split}
\end{equation*}
The Gronwall inequality implies
\begin{equation*}
 \mathbb{E} \sup_{s \leq t} |X_s^k - X_s^l|^2 \leq 4 \widehat{C} t \left|\frac{1}{k} - \frac{1}{l}\right| \exp \left\{ \left( 4C_1^2 S + 2 S + 8 K \right) t \right\} \,.
\end{equation*}
From this we can infer that there exists a process $(X_t)_{t \geq 0}$ such that
\begin{equation}\label{Xkconvergence}
 \mathbb{E} \sup_{s \leq t} |X_s - X_s^k|^2 \to 0 \text{ as } k \to \infty \,.
\end{equation}
It remains to be shown that $(X_t)_{t \geq 0}$ is indeed a solution to (\ref{multiplicativeSDE}). Observe that, by choosing a subsequence, we have $X_t^k \to X_t$ almost surely as $k \to \infty$ and thus, since $b$ is assumed to be continuous, we get
\begin{equation*}
 b(X_t^k - \frac{z}{k}) \to b(X_t) \text{ almost surely as } k \to \infty \,. 
\end{equation*}
But $b$ is bounded by the constant $\sqrt{M}$ and
\begin{equation*}
 \int_0^t \int_{\mathbb{R}^d} b(X_s^k - \frac{z}{k}) j(z) dz ds \leq \sqrt{M} \int_0^t \int_{\mathbb{R}^d} j(z) dz ds < \infty \,.
\end{equation*}
Therefore we get
\begin{equation*}
  \int_0^t \int_{\mathbb{R}^d} b(X_s^k - \frac{z}{k}) j(z) dz ds \to \int_0^t \int_{\mathbb{R}^d} b(X_s) j(z) dz ds = \int_0^t b(X_s) ds \text{ as } k \to \infty \text{ a.s.}
\end{equation*}
Moreover, using the It\^{o} isometry and (\ref{Xkconvergence}), we can easily prove that
\begin{equation*}
 \int_0^t \sigma (X_s^k) dW_s \to \int_0^t \sigma (X_s) dW_s 
\end{equation*}
and
\begin{equation*}
 \int_0^t \int_U g(X_{s-}^k,u) \widetilde{N}(ds,du) \to  \int_0^t \int_U g(X_{s-},u) \widetilde{N}(ds,du) 
\end{equation*}
almost surely (by choosing a subsequence), as $k \to \infty$, which finishes the proof.
\qed
\end{proof2}

Now we proceed with the proof of Lemma \ref{lemmaFeller}, which is needed to ensure existence of an invariant measure for the solution to (\ref{multiplicativeSDE}).

\begin{proof3}
 First observe that under our assumptions, we can use Lemma \ref{ymaEstimates} to get
 \begin{equation*}
  \mathbb{E} \sup_{s \leq t} |X_s|^2 \leq K_1(1+|x|^2) e^{K_2 t}
 \end{equation*}
for some constants $K_1$, $K_2 > 0$, where $(X_t)_{t \geq 0} = (X_t(x))_{t \geq 0}$ is a solution to (\ref{multiplicativeSDE}) with initial condition $x \in \mathbb{R}^d$. Hence, by the Chebyshev inequality, for any $\varepsilon > 0$ we can find $R > 0$ large enough so that for any $x \in \mathbb{R}^d$ with $|x| \leq R$ we have
\begin{equation}\label{lemma25proofSupremumBound}
 \mathbb{P} \left[ \sup_{s \leq t} |X_s(x)| \geq R \right] < \varepsilon \,.
\end{equation}
Now without loss of generality assume that $t \leq 1$ and fix $\varepsilon > 0$ and $R > 0$ like above. We can consider a solution $(X_t^R)_{t \geq 0}$ to the equation (\ref{multiplicativeSDE}) with the coefficients replaced by the truncated coefficients $b_R$, $\sigma_R$ and $g_R$ obtained from Lemma \ref{truncatingCoefficientsLemma} (note that the local boundedness assumption (\ref{localBoundedness}) in Lemma \ref{truncatingCoefficientsLemma} is satisfied due to the continuity of $b$ and the separate linear growth condition (\ref{linearGrowthSeparate}) for $\sigma$ and $g$, cf. the discussion at the beginning of this section). Then $b_R$, $\sigma_R$ and $g_R$ satisfy a global one-sided Lipschitz condition (\ref{GOLforTrunc}) with some constant $C(R) > 0$. Moreover, we have $X_s = X_s^R$ for $s \leq \tau_R$ with $\tau_R$ defined by 
\begin{equation}\label{defTau}
 \tau_R := \inf \{t > 0 : |X^R_t| \geq R \} \,.
\end{equation}
Thus for any $x$, $y \in \mathbb{R}^d$ with $|x| \leq R$ and $|y| \leq R$ and for any $\delta > 0$ we have
\begin{equation*}
 \mathbb{P}(|X_1(x) - X_1(y)| > \delta) \leq \varepsilon + \mathbb{P}(|X^R_1(x) - X^R_1(y)| > \delta) \leq \varepsilon + \frac{1}{\delta^2} \mathbb{E} |X^R_1(x) - X^R_1(y)|^2 \,,
\end{equation*}
where the first step follows from (\ref{lemma25proofSupremumBound}) and some straightforward calculations (see page 321 in \cite{albeverio} for details) and the second step is just the Chebyshev inequality. Now from the It\^{o} formula used similarly like in (\ref{ItoYma}) (cf. also the proof of Theorem \ref{globalOneSidedTheorem}, although here we need a~different local martingale than in the case where we estimate a supremum) we get
\begin{equation*}
\begin{split}
 |X^R_1(x) - X^R_1(y)|^2 &= |x-y|^2 + 2\int_0^1 \langle b_R(X^R_s(x)) - b_R(X^R_s(y)), X^R_s(x) - X^R_s(y) \rangle ds \\
 &+ 2\int_0^1 \langle \sigma_R(X^R_s(x)) - \sigma_R(X^R_s(y)), (X^R_s(x) - X^R_s(y)) dW_s \rangle \\
   &+ \int_0^1 \| \sigma_R(X^R_s(x)) - \sigma_R(X^R_s(y))\|_{HS}^2 ds \\
   &+2 \int_0^1 \int_U \Bigg\{ \langle g_R(X^R_{s-}(x),u) - g_R(X^R_{s-}(y),u), X^R_s(x) - X^R_s(y) \rangle \\
   &+ |g_R(X^R_{s-}(x),u) - g_R(X^R_{s-}(y),u)|^2 \Bigg\} \widetilde{N}(ds,du) \\
   &+ \int_0^1 \int_U  |g_R(X^R_{s-}(x),u) - g_R(X^R_{s-}(y),u)|^2 \nu(du) ds \\
   &\leq 2 C(R) \int_0^1 |X^R_s(x) - X^R_s(y)|^2 ds + M_t \,,
\end{split}
\end{equation*}
where we used the global one-sided Lipschitz condition (\ref{GOLforTrunc}) for $b_R$, $\sigma_R$ and $g_R$ and
\begin{equation*}
\begin{split}
 M_t &:= 2 \int_0^1 \int_U \Bigg\{ \langle g_R(X^R_{s-}(x),u) - g_R(X^R_{s-}(y),u), X^R_s(x) - X^R_s(y) \rangle  \\
 &+ |g_R(X^R_{s-}(x),u) - g_R(X^R_{s-}(y),u)|^2 \Bigg\} \widetilde{N}(ds,du) \\
 &+ 2\int_0^1 \langle \sigma_R(X^R_s(x)) - \sigma_R(X^R_s(y)), X^R_s(x) - X^R_s(y) dW_s \rangle 
 \end{split}
\end{equation*}
is a local martingale. Thus by a localization argument and the Gronwall inequality we get
\begin{equation}\label{lemma25proofExpectationBound}
 \mathbb{E} |X^R_1(x) - X^R_1(y)|^2 \leq A|x-y|^2 e^{B t}
\end{equation}
for some constants $A$, $B > 0$ and thus
\begin{equation}\label{lemma25proofProbabilityBound}
 \mathbb{P}(|X_1(x) - X_1(y)| > \delta) \leq \varepsilon + \frac{A}{\delta^2}|x-y|^2 e^{B t} \,.
\end{equation}
Once we have (\ref{lemma25proofProbabilityBound}), we proceed exactly like in the proof of Proposition 4.3 in \cite{albeverio}. Namely, we can show that for any sequence $x_n \to x$ in $\mathbb{R}^d$ we have $X_1(x_n) \to X_1(x)$ in probability. From this we infer that for any function $f \in \mathcal{C}_b(\mathbb{R}^d)$ we have
\begin{equation*}
 p_1f(x_n) \to p_1f(x) \,,
\end{equation*}
from which we get the desired Feller property of $(p_t)_{t \geq 0}$. Details of this last step can be found on page 321 in \cite{albeverio}.
 \qed
\end{proof3}

\begin{proof5}
 Note that (\ref{dissipativeGrowth}) obviously implies (\ref{GOLG}), hence under the assumptions of Theorem \ref{existenceInvariant}, the assumptions of Lemma \ref{lemmaFeller} are satisfied. Thus the semigroup $(p_t)_{t \geq 0}$ associated with the solution $(X_t)_{t \geq 0}$ of (\ref{multiplicativeSDE}) is Feller. Hence, if we can show (\ref{theorem45ineq}), then we can just use the Krylov-Bogoliubov method presented at the beginning of Section \ref{sectionInvariant} and conclude the proof.
In order to prove (\ref{theorem45ineq}), we apply the It\^{o} formula to $|X_t(x)|^2$ and then proceed like in the proof of Lemma \ref{lemmaFeller} presented above, where we apply the It\^{o} formula to obtain (\ref{lemma25proofExpectationBound}). However, unlike in (\ref{lemma25proofExpectationBound}), here we need to obtain the term $e^{Bt}$ with a negative constant $B$ in order to guarantee the boundedness in probability at infinity condition (\ref{boundedInProbability}). Thus we need to use the differential version of the Gronwall inequality and not the integral one (cf. Remark 2.3 in \cite{shao}). This is however not a problem, since by using (\ref{dissipativeGrowth}) and choosing a local martingale accordingly, we can obtain
\begin{equation*}
  \mathbb{E}|X_t(x)|^2 \leq \mathbb{E}|X_s(x)|^2 - 2K \int_s^t \left( \mathbb{E} |X_r(x)|^2 - \frac{M}{K} \right) dr 
\end{equation*}
for any $0 \leq s \leq t$. Thus by the differential version of the Gronwall inequality we have
\begin{equation*}
 \mathbb{E}|X_t(x)|^2 - \frac{M}{K} \leq |x|^2 e^{-2Kt} \,,
\end{equation*}
which gives (\ref{theorem45ineq}) and finishes the proof.
 \qed
\end{proof5}

\begin{proof4}
Under our assumptions we can prove that the coefficients of (\ref{multiplicativeSDE}) are locally bounded in the sense of (\ref{localBoundedness}), cf. the discussion at the beginning of Section \ref{sectionProofs}. Then combining Lemma \ref{truncatingCoefficientsLemma} and Theorem \ref{globalOneSidedTheorem}, we see that for every $R \geq 1$ there exists a unique strong solution $(X^R_t)_{t \geq 0}$ to the equation (\ref{multiplicativeSDE}) with the coefficients replaced by $b_R$, $\sigma_R$ and $g_R$ from Lemma \ref{truncatingCoefficientsLemma}. If we consider a sequence $\{ (X^n_t)_{t \geq 0} \}_{n \in \mathbb{N}}$ of such solutions and define stopping times $\{ \tau_n \}_{n \in \mathbb{N}}$ like in (\ref{defTau}), then we can show that for $n \leq m$ we have $\tau_n \leq \tau_m$ and consequently that 
\begin{equation*}
 \tau := \lim_{n \to \infty} \tau_n 
\end{equation*}
and
\begin{equation*}
 X_t := \lim_{n \to \infty} X_t^n \text{ almost surely on } [0,\tau]
\end{equation*}
are well defined. We just need to show that $(X_t)_{t \geq 0}$ is indeed a solution of (\ref{multiplicativeSDE}). However, by the construction of the coefficients in Lemma \ref{truncatingCoefficientsLemma} we can see that
\begin{equation*}
 X_{t \wedge \tau_n} = X^{n}_{t \wedge \tau_n} = X_0 + \int_0^{t \wedge \tau_n} b(X_s)ds + \int_0^{t \wedge \tau_n} \sigma(X_s)dW_s + \int_0^{t \wedge \tau_n} \int_{U} g(X_{s-},u) \widetilde{N}(ds,du) 
\end{equation*}
Therefore it remains to be shown that $\tau = \infty$ a.s., which can be done exactly as in Theorem 3.1 in \cite{albeverio}. Namely, using the global one-sided linear growth condition (\ref{GOLG}) we can show that $\mathbb{E}|X_{t \wedge \tau_n}|^2 \leq (\mathbb{E}|X_0|^2 + Kt) e^{Kt}$ for some constant $K > 0$ and then, after showing that $n^2 \mathbb{P}(\tau_n < t) \leq \mathbb{E}|X_{t \wedge \tau_n}|^2$, we see that $\tau_n \to \infty$ in probability and thus, via a subsequence, almost surely.
\qed
\end{proof4}

\section*{Acknowledgement}
I would like to thank Zdzis\l aw Brze\'{z}niak for discussions regarding the contents of \cite{albeverio}. This work was financially supported by the Bonn International Graduate School of Mathematics.

\end{document}